\documentclass[a4paper,12pt]{amsart}
\usepackage{amsmath,amsfonts,amssymb}
\usepackage{amsthm}
\usepackage{graphics,graphicx,epsf}
\usepackage[usenames]{color}
\usepackage{color}
\usepackage{hyperref} 
\usepackage{tikz}

\newtheorem{theorem}{Theorem}[section]
\newtheorem{definition}[theorem]{Definition}
\newtheorem{proposition}[theorem]{Proposition}
\newtheorem{lemma}[theorem]{Lemma}
\newtheorem{corollary}[theorem]{Corollary}
\newtheorem{remark}[theorem]{Remark}

\renewcommand{\proof}{{\noindent \bf Proof. \ }}
\newcommand{\sgn}{\textrm{sgn}}

\numberwithin{equation}{section}

\title[Nonlocal problems]{Upper semicontinuity for a class of nonlocal evolution equations with 
Neumann condition}

\author[F. D. M. Bezerra]{Flank  D. M. Bezerra}

\author[S. Sastre-Gomez]{Silvia Sastre-Gomez$^1$}\thanks{$^1$Partially supported by the projects MTM2012-31298 and Propesq Qualis A--2016.}

\author[S. H. Da Silva]{Severino H. da Silva$^2$}\thanks{$^2$Research partially
supported by CAPES/CNPq}

\date{\today}

\begin{document}

 \maketitle

\begin{abstract}
In this paper we consider the following nonlocal autonomous evolution equation in a bounded domain $\Omega$ in $\mathbb{R}^N$
\[
\partial_t u(x,t)    =- h(x)u(x,t)  + g \Big(\int_{\Omega} J(x,y)u(y,t)dy \Big) +f(x,u(x,t)) 
\]
where  $h\in W^{1,\infty}(\Omega)$, $g: \mathbb{R} \to \mathbb{R}$  and $f:\mathbb{R}^N\times\mathbb{R} \to \mathbb{R}$ are continuously differentiable function, and $J$ is a symmetric kernel; that is, $J(x,y)=J(y,x)$ for any $x,y\in\mathbb{R}^N$. Under additional suitable assumptions on $f$ and $g$, we study the asymptotic dynamics of the initial value problem associated to this equation in a suitable phase spaces. More precisely, we prove the existence,  and upper semicontinuity of compact global attractors with respect to kernel $J$.
\end{abstract}

\section{Introduction}

Nonlocal diffusion problems appear in many different areas such as neurology, 
ferromagnetism, medicine, biology and economics. Relevant models in Continuum Mechanics, Mathematical Physics, Biology and Economics are of nonlocal nature; for instance, Boltzmann equations in gas dynamics, Navier-Stokes equations in Fluid Mechanics, Keller-Segel model for Chemotaxis and Dynamic neural fields. There has been a deep 
study of existence, regularity of solutions and asymptotic dynamics of different 
nonlocal problems, c.f. \cite{Mazon_Rossi_2008, Rossi_libro, Berestycki, chasseigne2006asymptotic, da2009global, de1995travelling, de1994glauber, hutson, pereira2006global, ref6}. 

In this paper we are interested in the study of the  asymptotic dynamics of the  following nonlocal evolution problem with a nonlinear reaction term, 
in the sense of compact 
global attractors
\begin{equation}\label{prob}
\begin{cases}
\partial_t u(x,t) + h(x) u(x,t)  - g( K_J\,u(x,t))=f(x,u(x,t)), &  x \in  \Omega,\ t \geqslant  0, \\
u(x,0) = u_0 (x), & x \in  \Omega,
\end{cases}
\end{equation}
where $\Omega$ is a bounded  domain in $\mathbb{R}^N$, $N\geq1$; $h:\mathbb{R}^N\to \mathbb{R}$  is a continuously differentiable function with
\begin{equation}\label{Cond-h}
h(x)\geqslant h_0>0,\quad \mbox{and}\quad \partial_{x_i}h(x)\geqslant h_1>0,
\end{equation}
for any $i\in\{1,\ldots,N\}$ and $x\in\mathbb{R}^N$, and some constants $h_0, \,h_1$; $f:\mathbb{R}^N\times\mathbb{R} \to \mathbb{R}$  is  a continuously differentiable function that   satisfies the dissipative condition
\begin{equation}\label{dissip1}
|\xi(\cdot,s)|\leqslant k_f|s|+c_f
\end{equation}
where $\xi=f,\partial_2f$ or $\partial_1f$, and $k_f,\,c_f$  are strictly positive constants. 
 
The function $g: \mathbb{R} \to \mathbb{R}$ is  continuously differentiable  and satisfies the dissipative condition
\begin{equation}\label{dissip1g}
|\eta(s)|\leqslant k_g|s|+c_g,
\end{equation}
where $\eta=g$ or $g'$, with  $k_g,c_g$ strictly positive constants and 
\begin{equation}\label{hyp_h_0} 
k_f+k_g< h_0,
\end{equation}
this last condition will be used to prove the existence of compact attracting sets in $L^p(\Omega)$ for $1\leqslant p<\infty$. Furthermore,  $K_J$ is an integral operator 
\begin{equation}\label{mapK}
K_J \,v (x)  :=   \int_{\Omega} J(x,y)  v(y) d y.
\end{equation}
with symmetric kernel $J$; that is, $J(x,y)=J(y,x)$ for any $x,y\in \mathbb{R}^N$. We  assume without loss of generality that 
\begin{equation}\label{CondJ1}
\int_{\mathbb{R}^N} J(x,y) d  y = 1.
\end{equation}
Throughout this paper we will use the following notation, for  $1 \leqslant p 
\leqslant \infty$
\begin{equation}\label{CondJ2}  
\|J\|_{p}:= \sup_{ x \in \Omega } \|J(x,\cdot)\|_{L^{p} (\Omega) }< \infty.
\end{equation}

 The model \eqref{prob} includes models like the Ising spin system, which arises as one continuum limit of a probabilistic problem, where $u(x,t)$ denotes the magnetization density. Here, $f(x,u(x,t))$ describes the local rate of production or decrease of the magnetization density $u$ at the point $x$ at time $t \geq 0$. We are considering the case where $u(x,t)$ decays with speed $\frac{1}{h(x)}$ while it has a rate of production proportional  to a nonlinear function that depends on the points in a neighborhood of $x$ through the connectivity function $J$. The terms $h(x)$ and $f(x,u(x))$ are the main difference with respect to previous papers where the model  considered does not have the reaction term $f$ or they have a constant rate decay. 

The map $K_J$ given by \eqref{mapK} is well defined as a bounded linear operator in  several function spaces 
depending on the regularity assumed for $J$; for example, if $J$ satisfies 
\eqref{CondJ2} then $K_J$ is well defined in $L^p(\Omega)$ as shown below. 
We are interested in studying the asymptotic dynamics  of the problem 
\eqref{prob} in  $L^p(\Omega)$ with $1\leqslant p<\infty$, in the sense of 
compact global attractors.

The asymptotic behavior of solutions of  evolution equations with  nonlocal spatial terms has been extensively studied over the past years.  For instance, in  \cite{de1995travelling}, \cite{de2000critical}, \cite{de1994stability}, \cite{de1994uniqueness}, \cite{pereira2006global} and references therein, the authors study the following equation 
\begin{equation}\label{Prob1a}
\partial_t u(x,t) =- u(x,t)  + \tanh\big(\beta (J*u)(x,t)+ h\big),\  x \in  \mathbb{R},\ t \geqslant  0,
\end{equation}
with 
\[
(J*u) (x,\cdot)  :=   \int_{\mathbb{R}} J(x-y)u(y,\cdot)d y,
\]
where $\beta>1$, $J\in\mathcal{C}^1(\mathbb{R})$ is a non-negative even 
function with integral equal to 1 supported in $[-1,1]$, and 
$h$ is a positive constant. In \cite{de1994glauber}, it is shown that this 
equation arises as a continuum limit of one-dimensional Ising spin systems with 
Glauber dynamics and Kac potentials; $u$ represents a magnetization density 
and $\frac{1}{\beta}$ the temperature of the system. In \cite{pereira2006global}  
the author studied the existence  of global attractors and nonhomogeneous equilibria for  
\eqref{Prob1a}. 

In \cite{bezerra2012existence} and \cite{da2013finite} and references therein the authors consider the following Dirichlet nonlocal problem
\begin{equation}\label{Prob2a}
\partial_t u(x,t) =- u(x,t)  + g(\beta (J*u)(x,t)+\beta h),\  x \in\mathbb{R}^N,\ t \geqslant  0,
\end{equation}
where $u(x,t)=0$ for $x\in\mathbb{R}^N\backslash\Omega$ and $\Omega$ is a bounded smooth domain in $\mathbb{R}^N$, $N\geqslant1$, $\beta>1$, $h>0$, $g: \mathbb{R}  \to \mathbb{R}$ is a (sufficiently smooth) function and $J\in\mathcal{C}^1(\mathbb{R}^N)$ is a non-negative even function with integral equal to 1 supported in $\Omega$.  In \cite{bezerra2012existence} the authors studied the existence and characterization of global attractor for the Dirichlet problem \eqref{Prob2a} and in \cite{da2013finite} it has been proven the finite fractal dimensionality of this attractor. In this paper we study a variation of the equation \eqref{Prob2a}, considering the term $-h(x)u(x,t)$ instead of $-u(x,t)$ and adding the nonlinear reaction term $f(x,u(x,t))$.

In \cite{chasseigne2006asymptotic}, \cite{chasseigne2013nonlocal}, \cite{cortazar2005nonlocal},  \cite{cortazar2008approximate},  \cite{rodriguez2014nonlinear} the authors study the following  nonlinear nonlocal reaction-diffusion equation, and variations of it
\begin{equation}\label{Prob3a}
\partial_t u(x,t) =- u(x,t)  + (J*u)(x,t)+f(u(x,t)),\  x \in  \mathbb{R^N},\ t \geqslant  0, 
\end{equation}
where $J$ is a non-negative even function with integral equal to $1$ and $f: \mathbb{R}  \to \mathbb{R}$ is a (sufficiently smooth) function.

In \cite{amari1977dynamics}, \cite{da2012properties}, \cite{da2009global} and references therein the following equation associated to neural fields is considered
\begin{equation}\label{Prob4a}
\partial_t u(x,t) =- u(x,t)  + (J*(f\circ u))(x,t))+h,\  x \in  \mathbb{R},\ t \geqslant  0, 
\end{equation}
where $h>0$. This equation models the neuronal activity, and arise through a limiting argument from a discrete synaptically-coupled network of excitatory and inhibitory neurons. In (\ref{Prob4a}) the function $u(x,t)$ denotes the mean membrane potential of a patch of tissue located at position $x\in\mathbb{R}$ at time $t\geqslant 0$. The connection function $J(x-y)$ determines the coupling between the elements at position $x$ and position $y$. The function $f\circ u$ gives the neural firing rate, or average rate at which spikes are generated, corresponding to an activity level $u$. In [1] the author studied the  existence of a stationary travelling wave and in \cite{da2012properties} and \cite{da2009global}  the authors studied the existence of global attractor for  \eqref{Prob4a}. 

Motivated by these works  we consider in this paper a more general nonlocal 
equation  \eqref{prob}  defined in $L^p(\Omega)$,  with a nonlinear reaction term. In particular, for the problem \eqref{prob} we can consider the nonlocal Dirichlet 
problem and the nonlocal Neumann problem. Notice that the equation 
\eqref{prob} is defined only in $\Omega$ since for this equation we consider: 
the Dirichlet problem assuming that 
\begin{equation}\label{NCond}
u(x,t)= 0,\ \mbox{for all}\ x\in  \mathbb{R}^N \backslash \Omega  
\end{equation}
and the Neumann problem assuming that there is no flux across the 
boundary $\partial \Omega$. These two problems can be unified considering 
the equation \eqref{prob} defined in $\Omega$, (for more information see 
\cite{Rossi_libro, cortazar2008approximate, rodriguez2014nonlinear}). With this  we have a more general model in which we consider the Dirichlet and Neumann problem with just one model. From the viewpoint of mathematical analysis, we complement the analysis of the works \cite{bezerra2012existence,da2013finite,rodriguez2014nonlinear}, and \cite{rodriguez-bernal_sastre-gomez_2016} since in this paper we do not need the condition \eqref{NCond}. Moreover, we explore the Banach setting in our analysis of the nonlocal interactions in the PDE, and not only the Hilbert setting.

The paper is organized as follows. In Section \ref{wellposed} we prove the well posedness of the problem in $L^p(\Omega)$ with $1\leqslant p \leqslant\infty$. Moreover we show that the Cauchy problem (2.1) generates a nonlinear semigroup for $1\leqslant p \leqslant\infty$. In Section \ref{PullAttractors} we prove the existence  of  global attractor for $1\leqslant p <\infty$,  for this we show that there exists a compact attracting set in $L^p(\Omega)$. This result is proved under hypothesis \eqref{hyp_h_0} on $J$. Finally, in Section 4 we prove the upper semicontinuity of the global attractors  with respect to the kernel $J$.

\section{Well posedness of the problem}\label{wellposed}
In this section   we  show  that the problem (\ref{prob}) is well posed in $L^p(\Omega)$ for $1 \leqslant p \leqslant \infty$  in the sense described in  
Definition \ref{def_sol}. 

Let us rewrite problem \eqref{prob} as the Cauchy problem
\begin{equation}\label{CP}
\begin{cases}
 \partial_t u  = F(u),\ t>0, \\
 u(0)=u_{0}
\end{cases}
\end{equation}
where $F:L^p(\Omega) \to  L^p(\Omega)$ is the map defined by 
\begin{equation} \label{mapF}
F(u)=-h(\cdot)u+ g( K_J\,u)+f(\cdot,u).
\end{equation}
\begin{definition}\label{def_sol}

 A {\rm solution} of  \eqref{CP} in $[0,s)$ is a  continuous function $u:[0,s) \to 
 L^p(\Omega) $   such that $u(0)= u_{0} $,  the derivative with respect to $t$ exists and $\partial_t u(t,\cdot)$ belongs to $L^p(\Omega)$,
 and the differential equation in \eqref{CP} is satisfied for $ t \in [0, s)$.
\end{definition}

Below we give a result of some estimates for the nonlocal diffusion operator $K_J$ given by \eqref{mapK}, which  will be  used in the sequel.

\begin{lemma} \label{boundK}
    Let $K_J$ be the map defined  by (\ref{mapK}). 
  If $u \in   L^p{(\Omega)}, \ 1 \leqslant p \leqslant \infty$ then
 $K_J\,u \in L^{\infty}{(\Omega)}$, and we have
   \begin{equation}  \label{estimateLq}
 |K_J\, u (x)|  \leqslant   \|J\|_{p'}      \| u\|_{L^p{(\Omega)}}, \,
  \mathrm{for\ all}\  x \in \Omega,
 \end{equation}
  \textrm{ where $1\leqslant p' \leqslant \infty$
 is the conjugate exponent of $p$}. We also have 
  \begin{equation}  \label{estimateL1}
   \|K_J\, u\|_{L^{p}(\Omega)}  \leqslant  \|J\|_{1}
   \|  u\|_{L^{p}(\Omega)} \leqslant
   \|  u\|_{L^{p}(\Omega)}.
   \end{equation}
Moreover, if  $u \in   L^1{(\Omega)} $ then $K_J\,u \in L^{p}{(\Omega)}$, for  $1 \leqslant p \leqslant \infty$, and
\begin{equation}  \label{estimateLp}   
\|K_J\, u\|_{L^{p}(\Omega)}\leqslant  \|J\|_{p}\|  u\|_{L^{1}(\Omega)}.
\end{equation}
\end{lemma}

\proof
The estimate \eqref{estimateLq} follows easily from H\"{o}lder's inequality. The estimate \eqref{estimateL1} follows from
 the generalized  Young's inequality (see \cite{folland1995introduction}).
The proof of  (\ref{estimateLp}) is similar to
 (\ref{estimateL1}), but we include it here for the sake of completeness.
 Suppose $1< p < \infty$ and let $p'$ be its conjugate exponent.
 Then, by H\"{o}lder's inequality
 \begin{eqnarray*}
  |K_J\,u(x)| &\leqslant & \int_{\Omega} |J(x,y)| |u(y)|^{\frac{1 }{p}} |u(y)|^{\frac{1 }{p' }} 
dy      \\
  & \leqslant & \left(\int_{\Omega} |J(x,y)|^{p} |u(y)|   dy  \right)^{\frac{1}{p}}
       \left(\int_{\Omega}  |u(y)|  d y  \right)^{\frac{1}{p'}} \\
   & \leqslant& \| u \|_{L^{1}{(\Omega)}}^{\frac{1}{p'}}  \left(\int_{\Omega} |J(x,y)|^{p} |u(y)|   d y  \right)^{\frac{1}{p}}.  
\end{eqnarray*}
Raising both sides to the power $p$ and integrating in $\Omega$, we obtain
 \begin{eqnarray*}
  \int_{\Omega} |K_J\,u(x)|^p  d  x  & \leqslant &   \| u \|_{L^{1}{(\Omega)}}^{\frac{p}{p'}}
      \int_{\Omega}  \int_{\Omega} |J(x,y)|^{p} |u(y)|  d x
  d  y   \\
  & \leqslant &   \| u \|_{L^{1}{(\Omega)}}^{\frac{p}{p'}}  \int_{\Omega}   |u(y)|
  \int_{\Omega} |J(x,y)|^{p}   d x  d y  \\
  & \leqslant &   \| u \|_{L^{1}{(\Omega)}}^{\frac{p}{p'}}
  \| u \|_{L^{1}{(\Omega)}}   \|J  \|_p^p  \\
  & \leqslant &   \| u \|_{L^{1}{(\Omega)}}^{\frac{p+p'}{p'}}   \|J  \|_p^p.
  \end{eqnarray*}
  The inequality (\ref{estimateLp}) then follows
  by taking $p-th$ roots.

 The case $p=1$ is similar, and the case $p=\infty$ is trivial.
\qed

\begin{definition} \label{loclip}
If $E$ is a normed space, we say that a function  $F : E \to E$ is    \emph{locally Lipschitz  continuous  (or simply locally Lipschitz) } if,  for any $ x_0 \in   E$, there exists  a constant  $C$ and a neighborhood of $x_0$, $V = \{ x \in E \ : \   \|x-x_0  \|< b \}$, such that  if $ x$ and $ y$ belong to $V$   then
 $ \|F(x) - F(y) \| \leqslant  C \| x-y \| $; we say that $F$ is \emph{Lipschitz  continuous on bounded sets} if the neighborhood  $V$ in the previous definition can chosen as any bounded neighborhood in $E$.
 \end{definition}

\begin{remark}
  The two definitions in Definition \ref{loclip}  are equivalent if the normed space $ E  $ is locally compact.
\end{remark}

In the result below we prove that under suitable assumptions on $g$ and $f$, 
the map $F$ is Lipschitz continuous on bounded sets. 

\begin{proposition}\label{WellP}
Assume $g$  Lipschitz continuous on bounded sets, $f(\cdot,s)$ Lipschitz continuous on bounded sets on the second variable. Then, for each $1  \leqslant p \leqslant \infty$, the map $F$ defined by \eqref{mapF} is  Lipschitz continuous  on bounded sets.
\end{proposition}

\proof
Fix $u_0 \in L^p(\Omega)$. Let $V$ be the neighborhood of $u_0$ in $L^p(\Omega)$ given by,
\[
V: = \{ u \in  L^p(\Omega); \ \|u-u_0 \|_{L^{p}(\Omega)} < b \}.
\] 
From  (\ref{estimateLq}) in Lemma \ref{boundK}, it follows that   
\[
|  K_J\, u_0 (x) |  <    \|J\|_{p'} \|u_0 \|_{L^{p}(\Omega)},
\]
and for each $ u \in V$ and  $x \in \Omega$,
\[
\big| K_J\, u (x) - K_J\, u_0 (x) \big|  <    \|J\|_{p'}\, b.
\]
Let  $k_{V'}>0$ be  the Lipschitz constant of $g$ in the set
\[ 
V':=  \{  x \in   \mathbb{R}^N;  \  | x | \leqslant    \|J\|_{p'} (\|u_0 \|_{L^{p}(\Omega)} + b)  \}.
\]
If  $u, v \in V$,  then   for any $x \in \Omega$
\[
 |g( K_J\,u(x) )  - g(K_J\,v(x) ) | \leqslant k_{V'}   | K_J\,u(x)  - K_J\,v(x)  | .
 \]
For   $ 1 \leqslant p < \infty$, it follows then, from (\ref{estimateL1}) in Lemma \ref{boundK}, we have
\[
 \begin{split}
 \|g( K_J\,u)  - g( K_J\,v )\|_{L^p(\Omega)}
  & \leqslant    \left(
 \int_{\Omega}    k_{V' }  \big| K_J\,u(x)  - K_J\,v(x)  \big|^p   dx\right)^{\frac{1}{p}} \\
 &  =  
 k_{V'}  \left(  \int_{\Omega}     \big| K_J\, (u-v)(x) \big|^p   dx\right)^{\frac{1}{p}}  \\
  &  =  
 k_{V'}   \| K_J\, (u-v) \|_{L^p(\Omega)} \\
 &  \leqslant   k_{V'}   \| u-v \|_{L^p(\Omega)}.
\end{split}
\]

If $p =\infty$ the same inequality follows immediately from \eqref{estimateL1} in Lemma \ref{boundK}. Thus,  the map
\[
u \in L^p(\Omega) \mapsto  g(K_J\, u) \in L^p(\Omega) 
\] 
 is  Lipschitz continuous  on the set $V$. We also have that for any $u, v \in V$ and $x\in\Omega$,
\[
 \begin{split}
 \|f(\cdot, u)  - f(\cdot,v )\|_{L^p(\Omega)} & =\left( \int_{\Omega} \big| f(x,u(x))-f(x,v(x))\big|^p   dx\right)^{1/p} \ \\
  & \leqslant    \left(
 \int_{\Omega}    k_{V }^p  \big| u(x)  - v(x)  \big|^p   dx\right)^{1/p} \\
 &  =   k_{V}   \| u-v \|_{L^p(\Omega)},
\end{split}
\]
 and, we obtain that for any $u, v \in V$,
\[
\| h(u -v) \|_{L^p(\Omega)} \leqslant \|h\|_{L^{\infty}(\Omega)}\|u -v\|_{L^p(\Omega)},
\]
and so
\[
\begin{split}
\|  F(u)-F(v) \|_{L^p(\Omega)}&\leqslant \|  h (u-v)  \|_{L^p(\Omega)} +\| g(K_J\,u)-g(K_J\,v)  \|_{L^p(\Omega)}+\| f(u)-f(v) \|_{L^p(\Omega)}\\
&\leqslant C\| u-v  \|_{L^p(\Omega)},
\end{split}
\]
where $C=C\big(V,V',\|h\|_{L^{\infty}(\Omega)}\big)>0$ is constant, and therefore, the map $F$ is  Lipschitz continuous on the bounded set $V$.
\qed
\\

From the result above, it follows from well known results that  the  problem  \eqref{CP} has a local solution for any initial data in $L^p(\Omega)$. 
 For the  global existence, we need the following result (c.f. \cite[Theorem 5.6.1]{ladas1972differential}).

\begin{theorem}\label{ladas} 
Let $X$ be a Banach space, and assume $f: [ t_0, +\infty)  \times X \to X$  is continuous and satisfies
\[
 \|f(t,u) \|  \leqslant \Phi\big(t, \|u  \|\big),\; \mbox{ for all }\; (t,u) \in    [t_0, +\infty) \times X,
\]
where $\Phi: [t_0, +\infty) \times [0,+\infty)  \to  [0,+\infty) $ is continuous, 
$ \Phi(t,r)$  is  non decreasing in $r \geqslant 0$, for each $t \in [ t_0, +\infty) $ and the maximal solution $r(t;t_0,r_0) $ of the scalar initial value problem
 \[
 \begin{cases}  
 \dfrac{dr}{dt}= \Phi(t,r),\ t \geqslant t_0,\\
 r(t_0) = r_ 0, 
\end{cases}  
 \]
 exists for all  $t\in [t_0, +\infty)$. Then the maximal interval of existence of any  solution
 $u(t;t_0,y_0)$ of the initial value problem
 \[ 
 \begin{cases}
 \dfrac{du}{dt} = f(t,u),\ t \geqslant t_0,\\
 u(t_0) = u_0,
\end{cases} 
 \]
 is given by $[ t_0, +\infty)$.

\end{theorem}

\begin{corollary} \label{globalexist}
Under hypotheses in Proposition \ref{WellP}  the   problem  (\ref{CP})  has a unique global solution for any initial
 condition in  $L^p(\Omega)$, for $1\le p\le \infty$, which is  given by
 \begin{equation}\label{EP_1}
u(x,t) =
 e^{-th(x)}u_0(x)+\displaystyle\int_0^t e^{-(t-s)h(x)}\big[g(K_J\,u(x,s))+f(x,u(x,s))\big]ds. 
  \end{equation}

\end{corollary}

\proof
From Proposition \ref{WellP}, it follows that the right-hand-side of \eqref{CP} is  Lipschitz continuous on bounded sets of
 $L^p(\Omega)$ and, therefore, the Cauchy problem \eqref{CP} is well posed in $ L^p(\Omega)$, with a unique local solution $u(t,x)$,
 given by \eqref{EP_1}, (c.f. \cite{daleckii2002stability}). 
 

Thanks to \eqref{estimateL1} in  Lemma \ref{boundK}, \eqref{dissip1} and \eqref{dissip1g}, if $ 1 \leqslant p < \infty$,  we obtain  the following estimate
\[
\|f(\cdot,u)\|_{L^p(\Omega)}   \leqslant c_f |\Omega|^{\frac{1}{p}} + k_f \|   u \|_{L^p(\Omega)}.
\]
and
\[
\|g(K_J\,u)\|_{L^p(\Omega)}   \leqslant c_g |\Omega|^{\frac{1}{p}} + k_g \|   u \|_{L^p(\Omega)}.
\]

For $p= \infty$, using similar arguments (or by making $p \to \infty$), we obtain 
that
\[
  \| f(\cdot,u )\|_{L^{\infty}(\Omega)} \leqslant   c_f  + k_f \|    u  \|_{L^{\infty}(\Omega)}.
\]
and
\[
  \| g(K_J\,u )\|_{L^{\infty}(\Omega)} \leqslant   c_g  + k_g \|    u  \|_{L^{\infty}(\Omega)}.
\]

Now defining the function $\Phi: [t_0, +\infty) \times [0,+\infty)  \to  [0,+\infty)$ in Theorem \ref{ladas} by
\[
\Phi(t,r)= (c_f+c_g)|\Omega|^{\frac{1}{p}}+(\|h\|_{L^{\infty}(\Omega)}+k_f+k_g) r,
\]
for any $(t,r)\in  [t_0, +\infty) \times [0,+\infty)$, it follows that the problem \eqref{CP} satisfies the hypotheses of  Theorem \ref{ladas}  and the global existence  follows immediately. The variation of constants formula can be verified by direct differentiation. \qed
\\

The result below will be used to prove that the map $F$ defined by 
\eqref{mapF} is continuous Fr\'echet differentiable. The proof can be seen in 
\cite{rall2014nonlinear}.

\begin{proposition}\label{Prop-Rall}
Let $Y$ and $Z$ be normed linear spaces, $F:Y\to Z$ a map and suppose that the Gateaux derivative of $F$, $DF:Y\to\mathcal{L}(Y,Z)$ exists and is continuous at $y\in Y$. Then the Frech\'et derivative $F' $ of $F$ exists and is continuous at $y$.
\end{proposition}

\begin{proposition} \label{C1flow}
The map $F$ defined by 
\eqref{mapF} is continuously Frech\'et differentiable with derivative $DF: L^p(\Omega)\to  \mathcal{L}(L^{p}(\Omega), L^{1}(\Omega))$ given by
\begin{equation}\label{DiffF}
DF(u)v(x):=h(x)v(x)+ g'( K_J\,u(x))(K_J\,v(x))+\partial_2f(x,u(x))v(x),
\end{equation}
for any $u,v\in L^p(\Omega)$, $1\le p\le \infty$ and $x\in\Omega$.
\end{proposition}

\proof From a simple computation, using that $f$ and $g$ are continuously differentiable, it follows that the Gateaux's derivative of $F$ is given by
\[
DF(u)v(x):=h(x)v(x)+ g'( K_J\,u(x))(K_J\,v(x))+\partial_2f(x,u(x))v(x), 
\]
for any $u,v\in L^p(\Omega)$ and $x\in\Omega$.

To apply Proposition \ref{Prop-Rall} we will consider $Y=L^p(\Omega)$ and $Z=L^1(\Omega)$. Since $F:L^p(\Omega)\to L^p(\Omega)$ we have that  $F:L^p(\Omega)\to L^1(\Omega)$. Thus, the first hypothesis in Proposition \ref{Prop-Rall} is satisfied. The operator $DF(u)$ is clearly  a linear operator  from $L^p(\Omega)$ to $L^1(\Omega)$. Suppose $1 \leqslant p < \infty $, then for $u \in L^p(\Omega)$,  we have
\begin{equation}\label{In_0}
\|hv\!+ \! g'( K_J\,u)(K_J\,v)\!+ \!\partial_2 f(\cdot,u)v \|_{ L^1(\Omega) }\! \leqslant\! \|hv\|_{ L^1(\Omega) }\!+ \!\|g'( K_J\,u)(K_J\,v)\|_{ L^1(\Omega) }\!+\!\|\partial_2f(\cdot,u)v \|_{ L^1(\Omega) },
\end{equation}
where
\begin{equation}\label{In_1}
\|hv\|_{ L^1(\Omega) }\leqslant |\Omega|^{1\over p'}\|h\|_{L^{\infty}(\Omega)}\|v\|_{ L^p(\Omega) }.
\end{equation}
From \eqref{dissip1g}, Lemma \ref{boundK} and Minkowski's inequality, we obtain that
\[
\begin{split}
\|g'( K_J\,u)(K_J\,v)\|_{ L^1(\Omega) }&= \int_\Omega\big|g'( K_J\,u(x))\big|\big|K_J\,v(x)\big|dx\\
&\leqslant \left(\int_\Omega \big|g'( K_J\,u(x))\big|^{p'}dx\right)^{\frac{1}{p'}} \left(\int_\Omega  \big|K_J\,v(x)\big|^{p}dx\right)^{\frac{1}{p}}\\
&\leqslant \left( \int_\Omega \big|k_g| K_J\,u(x)|+c_g\big|^{p'}dx\right)^{\frac{1}{p'}} \|v\|_{L^p(\Omega)}\\
&\leqslant  \big(k_g\|u\|_{L^p(\Omega)} +|\Omega|^{\frac{1}{p'}}c_g\big) \|v\|_{L^p(\Omega)}.
\end{split}
\]

From \eqref{dissip1}  we have 
\begin{equation}\label{In_3}
\begin{split}
\|\partial_2f(\cdot,u)v \|_{ L^1(\Omega) }&= \int_\Omega|\partial_2f(x,u(x))v(x)|dx\\
&\leqslant \left(\int_\Omega \big|\partial_2f(x,u(x))\big|^{p'}dx\right)^{\frac{1}{p'}} \left(\int_\Omega  |v(x)|^{p}dx\right)^{\frac{1}{p}} \\
&\leqslant \left(\int_\Omega \big|k_f|u(x)|+c_f\big|^{p'}dx\right)^{\frac{1}{p'}} \left(\int_\Omega  |v(x)|^{p}dx\right)^{\frac{1}{p}} \\
&\leqslant  \left(k_f\|u\|_{L^p(\Omega)}+c_f|\Omega|^{\frac{1}{p'}}\right) \|v\|_{L^p(\Omega)},
\end{split}
\end{equation}
where $p'$ is the conjugate exponent of $p$.

Hence using \eqref{In_0} jointly with \eqref{In_1} and \eqref{In_3} we conclude that 
\[
\big\|hv+ g'( K_J\,u)(K_J\,v)+\partial_2f(\cdot,u)v \big\|_{ L^1(\Omega) } \leqslant C\|v\|_{ L^p(\Omega) }
\]
where $C=C(\|h\|_{L^{\infty}(\Omega)},\|u\|_{L^p(\Omega)},|\Omega|,g,f,p)>0$, showing that $DF(u)$ is a bounded linear map from $L^p(\Omega)$ to $L^1(\Omega)$.  Suppose now that $u_n$, $u$ and $v$  belong to $L^p(\Omega)$, $ 1 \leqslant p  < \infty$ with $ \|  u_n - u  \|_{L^p(\Omega)} \to 0 $ as $n\to\infty$. From Lemma \ref{boundK}  and H\"{o}lder's inequality, it follows  that
\[
\begin{split}
&\| DF(u_n)v - DF(u) v\|_{L^1(\Omega)} \\
& \leqslant \|g'( K_J\,u_n)(K_J\,v) -g'( K_J\,u)(K_J\,v) \|_{L^1(\Omega)}+\| \partial_2f(\cdot,u_n)v-\partial_2f(\cdot,u)v \|_{L^1(\Omega)}\\
& \leqslant \|v\|_{L^p(\Omega)}\|g'( K_J\,u_n)-g'( K_J\,u)\|_{L^1(\Omega)}+\| \partial_2f(\cdot,u_n)v-\partial_2f(\cdot,u)v \|_{L^1(\Omega)}\\
& \leqslant \|v\|_{L^p(\Omega)}\|g'( K_J\,u_n)-g'( K_J\,u)\|_{L^1(\Omega)}+\| \partial_2f(\cdot,u_n)-\partial_2f(\cdot,u) \|_{L^{p'}(\Omega)}\|v\|_{L^{p}(\Omega)}.
\end{split}
\]    

Thus to prove continuity of the derivative, we note that $ \|  u_n - u  \|_{L^p(\Omega)} \to 0 $ as $n\to\infty$ and we see that $ u_n \to u $ a.e in $\Omega$ as $n\to\infty$, and this implies $Ku_n \to Ku $ in $\Omega$ as $n\to\infty$, and therefore, there exists a bounded set $B\subset\mathbb{R}$ such that $g'$ is Lipschitz continuous on $B$ and
\[
\begin{split}
\|g'( K_J\,u_n)-g'( K_J\,u)\|_{L^1(\Omega)} &=\int_\Omega |g'( K_J\,u_n)-g'( K_J\,u)|dx\\
&\leqslant C_B\int_\Omega |K_J\,(u_n-u) |dx\\
&\leqslant C\|u_n-u\|_{L^p(\Omega)},
\end{split}
\] 
where $C>0$ depends on the Lipschitz constant of $g'$ on $B$. Moreover, since $ u_n \to u $ a.e in $\Omega$ as $n\to\infty$, there exists a bounded set $D\subset\mathbb{R}$ such that $\partial_2f(x,\cdot)$ is Lipschitz continuous on $D$ and
\[
\begin{split}
\| \partial_2f(\cdot,u_n)-\partial_2f(\cdot,u) \|_{L^{p'}(\Omega)}&=\left(\int_\Omega| \partial_2f(\cdot,u_n)-\partial_2f(\cdot,u)|^{p'} dx\right)^{\frac{1}{p'}}\\
&\leqslant C_D\left(\int_\Omega |u_n-u|^{p'} dx\right)^{\frac{1}{p'}}\\
&\leqslant C_D\|u_n-u\|_{L^{p}(\Omega)}.
\end{split}
\]
The previous computations are analogous for the case $p=\infty$. Therefore, it  follows from  Proposition \ref{Prop-Rall} above that $F$ is Fr\'echet differentiable  with  continuous  derivative. \qed


From the results above we have the following result.

\begin{theorem}
For each  $u_0 \in L^p(\Omega)$,   for $1\le p\le \infty$,  the unique solution of \eqref{CP} with initial condition $u_0$ exists for all $t\geqslant0$. Namely, this solution $(x,t)\mapsto u(x,t)$,  defined by \eqref{EP_1}, gives rise to a family of nonlinear $C^0$-semigroup in $ L^p(\Omega)$, $\{S_J(t);t\geqslant0\}$, which is given by
\[
S_J(t)u_0(x):=u(x,t),
\]
for any $x\in\Omega$ and $t\geqslant0$.
\end{theorem}
 
The notation $S_J(\cdot)$ refers to dependence of the semigroups in relation to the kernel $J$, this dependence will be explored in the following sections.

\section{Existence of global attractor}\label{PullAttractors}

In this section we prove the existence of global attractor $\mathcal{A}_J$ in $L^p(\Omega)$ for the nonlinear semigroup $ \{S_J(t); t \geqslant 0\} $ for $ 1 \leqslant p  <  \infty $, using \cite[Theorem 2.1]{carvalho2012attractors}. More precisely ,we will prove that the semigroup has a compact attracting set in $L^p(\Omega)$ with $ 1 \leqslant p  <  \infty $.

\begin{lemma}\label{L_PullAbs}
Suppose that the hypotheses \eqref{Cond-h}, \eqref{dissip1} and \eqref{dissip1g} hold and the constants $h_0$, $k_f$ and $k_g$ satisfy \eqref{hyp_h_0}.  Then the ball of $L^p(\Omega), \ 1 \leqslant p <  \infty$, centered
at the origin with  radius $r_\delta$ defined by
\begin{equation}\label{DefR0}
r_\delta=\dfrac{1}{h_0 -k_f-k_g } (c_f+c_g)(1+\delta)\max\{1,|\Omega|\},
\end{equation}
which we denote by $\mathcal{B}(0;r_\delta)$, where $c_f$ and  $c_g$ are the constants in \eqref{dissip1} and \eqref{dissip1g}, respectively, and $\delta$ is any positive constant,   absorbs bounded subsets of $L^p(\Omega)$  with respect to the nonlinear semigroup $S_J(\cdot)$ generated by \eqref{CP}.
\end{lemma}

\proof  
Let $u(x,t)$ be the solution of \eqref{CP} with initial condition $u_0\in L^p(\Omega)$ then, for $ 1 \leqslant p <\infty$, we have
 \begin{equation}\label{mkq1}
\begin{split}
\frac{d}{dt}\int_{\Omega}|u(x,t)|^{p}dx & =  \int_{\Omega} p |u(x,t)|^{p-1} \sgn(u(x,t))\,\partial_tu(x,t) dx  \\
&= -p\int_{\Omega}h(x)|u(x,t)|^{p}dx +p\int_{\Omega}|u(x,t)|^{p-1} \sgn(u(x,t))\,g( K_J\,u(x,t)) dx\\
&+p\int_{\Omega}|u(x,t)|^{p-1} \sgn(u(x,t))\,f(x,u(x,t)) dx.
\end{split}
\end{equation}

Note that from \eqref{Cond-h}, we have
\begin{equation}
p\,h_0 \int_{\Omega}|u(x,t)|^{p}dx\leqslant p\int_{\Omega}h(x)|u(x,t)|^{p}dx.
\end{equation}

Using H\"{o}lder's inequality, \eqref{estimateL1} in Lemma \ref{boundK}  and  \eqref{dissip1g},  we obtain that
 \begin{equation}\label{mkq2}
\begin{split}
&\int_{\Omega}|u(x,t)|^{p-1} \sgn(u(x,t))g( K_J\,u(x,t)) dx\\
 & \leqslant  \left(\int_{\Omega} |u(x,t)|^{p'(p-1)} dx\right)^{\frac{1}{p'}} \left(\int_{\Omega} | g( K_J\,u(x,t))|^pdx \right)^{\frac{1}{p}}  \\
 & \leqslant  \left(\int_{\Omega} |u(x,t)|^{p}dx\right)^{\frac{1}{p'}} \left(\int_\Omega \big|k_g| K_J\,u(x,t)|+c_g\big|^pdx\right)^{\frac{1}{p}}  \\
 &\leqslant  \|u(\cdot,t)\|_{L^{p}(\Omega)}^{p-1} \left(k_g \|u(\cdot,t)\|_{L^{p}(\Omega)} + c_g|\Omega|^{\frac{1}{p}}\right),
\end{split}
\end{equation}
where $p'$ is the conjugate exponent of $p$. 
We also have
\begin{equation}\label{dp3d93}
\begin{split}
&\int_{\Omega}|u(x,t)|^{p-1} \sgn(u(x,t))f(x,u(x,t)) dx\\
&\leqslant \left(\int_{\Omega} |u(x,t)|^{p'(p-1)} dx\right)^{\frac{1}{p'}} \left(\int_{\Omega} | f(x,u(x,t))|^pdx \right)^{\frac{1}{p}}  \\
& \leqslant  \left(\int_{\Omega} |u(x,t)|^{p}dx\right)^{\frac{1}{p'}} \left(\int_\Omega \big|k_f| u(x,t)|+c_f\big|^pdx\right)^{\frac{1}{p}}  \\
&\leqslant  \|u(\cdot,t)\|_{L^{p}(\Omega)}^{p-1} \left(k_f \|u(\cdot,t)\|_{L^{p}(\Omega)} + c_f|\Omega|^{\frac{1}{p}}\right).
\end{split}
\end{equation}

Hence, from \eqref{mkq1}-\eqref{dp3d93} we conclude that
\[
\begin{split}
&\frac{d}{dt}\|u(\cdot,t)\|_{L^{p}(\Omega)}^{p}\\
 &\leqslant -ph_0\|u(\cdot,t)\|_{L^{p}(\Omega)}^{p}+p(k_f+k_g)\|u(\cdot,t)\|_{L^{p}(\Omega)}^{p} +p(c_f+c_g){|\Omega|}^{\frac{1}{p}} \|u(\cdot,t)\|^{p-1}_{L^{p}(\Omega)}\\
&= p\|u(\cdot,t)\|_{L^{p}(\Omega)}^{p} \left[-h_0 +k_f+k_g  +\frac{|\Omega|^{\frac{1}{p}}}{\|u(\cdot,t)\|_{L^{p}(\Omega)}}(c_f+c_g)\right]\\
&\leqslant p\|u(\cdot,t)\|_{L^{p}(\Omega)}^{p} \left[-h_0 +k_f+k_g  +\frac{\max\{1,|\Omega|\}}{\|u(\cdot,t)\|_{L^{p}(\Omega)}}(c_f+c_g)\right].
\end{split}
\]

Let   $\varepsilon =h_0 -(k_f+k_g)  >0$. Then, if we consider
\[
\|u(\cdot,t)\|_{L^{p}(\Omega)} \geqslant \dfrac{1}{\varepsilon}  (c_f+c_g)(1+\delta)\max\{1,|\Omega|\},
 \]
 we have
 \[
\begin{split}
\frac{d}{dt}\|u(\cdot,t)\|_{L^{p}(\Omega)}^{p} &\leqslant
p \|u(t,\cdot)\|_{L^{p}(\Omega)}^{p}\left(-\varepsilon +
\frac{\varepsilon}{1+\delta}\right)\\
&=- \frac{\delta \varepsilon p}{1+\delta} \|u(t,\cdot)\|_{L^{p}(\Omega)}^{p}.
\end{split}
\]

Therefore
\begin{equation}\label{boundsol}
\begin{split}
\|u(\cdot,t)\|_{L^{p}(\Omega)}^{p} &\leqslant e^{-\frac{ \delta \varepsilon p}{(1+\delta)}t} \| u_0\|_{L^{p}(\Omega)}^{p} \\
&=e^{- \frac{\delta p}{(1+\delta)}(h_0 -k_f-k_g)t}\| u_0\|_{L^{p}(\Omega)}^{p}.
\end{split}
\end{equation}
From this, the  result follows easily for $1 \leqslant p < \infty$, and this completes the proof of the lemma.
\qed

\begin{remark}
From \eqref{boundsol} it follows that   the dissipativity of the semigroup in the space $L^{p}(\Omega)$ with $1 \leqslant p < \infty$ actually is independent of $p$, and therefore, we can assure the dissipativity of the semigroup in the space $L^{\infty}(\Omega)$, and this allows us to consider $u(x,t)$ essentially bounded in $\Omega$ for $t$ sufficiently large.
\end{remark}

\begin{theorem}\label{TExistGloAt}
In addition to the  conditions of  Lemma \ref{L_PullAbs},  suppose that for each $i\in\{1,\ldots,N\}$ and for each $ 1 \leqslant r \leqslant \infty$, we have
\[
\|\partial_{x_i} J\|_{r}:= \sup_{ x \in \Omega } \|\partial_{x_i}J(x,\cdot)\|_{L^{r} (\Omega) }< \infty,
\]  
and if $h_0>k_fr_\delta+c_f$, where $r_\delta$ is given by \eqref{DefR0}, then there exists a global attractor $\mathcal{A}_J$ for the nonlinear semigroup $\{S_J(t); t\geqslant 0\}$ generated by \eqref{CP} in $L^p(\Omega)$,   for $1\le p< \infty$.
\end{theorem}

\proof  
We will prove that the semigroup has a compact attracting set in $L^p(\Omega)$ with $ 1 \leqslant p  <  \infty $. If  $u(x,t)$ is the solution of \eqref{CP} with initial condition $u_0$ then, for $ 1 \leqslant p <\infty$ and $t$ sufficiently large we have
 \begin{equation}\label{mkqa2da1}
\begin{split}
&\frac{d}{dt}\int_{\Omega}|\partial_{x_i}u(x,t)|^{p}dx\\
 & =  p\int_{\Omega}  |\partial_{x_i} u(x,t)|^{p-1} \sgn(\partial_{x_i}u(x,t))\partial_{x_i}\partial_tu(x,t) dx  \\
&= -p\int_{\Omega}\partial_{x_i}h(x)|\partial_{x_i}u(x,t)|^{p-1}\sgn(\partial_{x_i}u(x,t))u(x,t)dx -p\int_{\Omega} h(x)|\partial_{x_i} u(x,t)|^{p} dx\\
& +p\int_{\Omega}|\partial_{x_i}u(x,t)|^{p-1} \sgn(\partial_{x_i}u(x,t))g'( K_J\,u(x,t)) ( K_{\partial_{x_i}J}u)(x,t) dx\\
&+p\int_{\Omega}|\partial_{x_i}u(x,t)|^{p-1} \sgn(\partial_{x_i}u(x,t))\partial_1f(x, u(x,t)) dx\\
&+p\int_{\Omega}|\partial_{x_i}u(x,t)|^{p-1} \sgn(\partial_{x_i}u(x,t))\partial_2f(x, u(x,t)) \partial_{x_i} u(x,t) dx.
\end{split}
\end{equation}


Now note that by H\"older's inequality and \eqref{Cond-h} we obtain that
\begin{equation}\label{ad33fn9}
\begin{split}
&-p\int_{\Omega}\partial_{x_i}h(x)|\partial_{x_i}u(x,t)|^{p-1}\sgn(\partial_{x_i}u(x,t))u(x,t)dx\\
&\leqslant p\|h\|_{W^{1,\infty}(\Omega)}\|\partial_{x_i}u(\cdot,t)\|_{L^{p}(\Omega)}^{p-1}\|u(\cdot,t)\|_{L^{p}(\Omega)}\\
&\leqslant pr_\delta\|h\|_{W^{1,\infty}(\Omega)}\|\partial_{x_i}u(\cdot,t)\|_{L^{p}(\Omega)}^{p-1}.
\end{split}
\end{equation}

Using H\"{o}lder's inequality and the estimate \eqref{estimateL1} in Lemma \ref{boundK} we get 
\[
\begin{split}
&\int_{\Omega}|\partial_{x_i}u(x,t)|^{p-1} \sgn(\partial_{x_i}u(x,t))g'( K_J\,u(x,t))( K_{\partial_{x_i}J}u)(x,t) dx\\
 & \leqslant  \left(\int_{\Omega} |\partial_{x_i}u(x,t)|^{p'(p-1)} dx\right)^{\frac{1}{p'}} \left(\int_{\Omega} | g'( K_J\,u(x,t))|^p|  ( K_{\partial_{x_i}J}u)(x,t) |^pdx \right)^{\frac{1}{p}}  \\
  & \leqslant \|u(\cdot,t)\|_{L^{p}(\Omega)}  \|\partial_{x_i} J\|_{p'}\left(\int_{\Omega} |\partial_{x_i}u(x,t)|^{p'(p-1)} dx\right)^{\frac{1}{p'}} \left(\int_{\Omega} | g'( K_J\,u(x,t))|^p  dx \right)^{\frac{1}{p}},
\end{split}
\]
where $p'$ is the conjugate exponent of $p$, and by the dissipative condition \eqref{dissip1},  we have that for $t$ sufficiently large
 \begin{equation}\label{mkqsdasd2}
\begin{split}
&\int_{\Omega}|\partial_{x_i}u(x,t)|^{p-1} \sgn(\partial_{x_i}u(x,t))g'( K_J\,u(x,t))( K_{\partial_{x_i}J}u)(x,t) dx\\
 & \leqslant  \|u(\cdot,t)\|_{L^{p}(\Omega)} \|\partial_{x_i} J\|_{p'} \left(\int_{\Omega} |\partial_{x_i}u(x,t)|^{p}dx\right)^{\frac{1}{p'}} \left(\int_\Omega |k_g| K_J\,u(x,t)|+c_g|^pdx\right)^{\frac{1}{p}}  \\
 &\leqslant  \|\partial_{x_i} J\|_{p'}\|u(\cdot,t)\|_{L^{p}(\Omega)} \|\partial_{x_i}u(\cdot,t)\|_{L^{p}(\Omega)}^{p-1} \left(k_g \|u(\cdot,t)\|_{L^{p}(\Omega)} + c_g|\Omega|^{\frac{1}{p}}\right)\\
 &\leqslant r_\delta \|\partial_{x_i} J\|_{p'}\big(k_g r_\delta + c_g|\Omega|^{\frac{1}{p}} \big)\|\partial_{x_i}u(\cdot,t)\|_{L^{p}(\Omega)}^{p-1} .
\end{split}
\end{equation}
We also have that
\begin{equation}\label{asda23s}
\begin{split}
&\int_{\Omega}|\partial_{x_i}u(x,t)|^{p-1} \sgn(\partial_{x_i}u(x,t))\partial_1f(x, u(x,t)) dx\\
&\leqslant    \left(\int_{\Omega} |\partial_{x_i}u(x,t)|^{p'(p-1)} dx\right)^{\frac{1}{p'}} \left(\int_{\Omega} | \partial_1f(x, u(x,t))|^p  dx \right)^{\frac{1}{p}}  \\
&\leqslant    \left(\int_{\Omega} |\partial_{x_i}u(x,t)|^{p} dx\right)^{\frac{1}{p'}} \left(\int_{\Omega} | k_f|u(x,t)|+c_f|^p  dx \right)^{\frac{1}{p}}  \\
&\leqslant  \big(k_f r_\delta+ c_f|\Omega|^{\frac{1}{p}}\big)  \|\partial_{x_i}u(\cdot,t)\|_{L^{p}(\Omega)}^{p-1}.
\end{split}
\end{equation}
We also have
\[
\begin{split}
&\int_{\Omega}|\partial_{x_i}u(x,t)|^{p-1} \sgn(\partial_{x_i}u(x,t))\partial_2f(x, u(x,t)) \partial_{x_i} u(x,t) dx\\
&=\int_{\Omega}|\partial_{x_i}u(x,t)|^{p} \sgn(\partial_{x_i}u(x,t))\partial_2f(x, u(x,t)) dx\\
& \leqslant \int_{\Omega}|\partial_{x_i}u(x,t)|^{p}|k_f|u(x,s)|+c_f|dx,
\end{split}
\]
and so
\begin{equation}\label{dp3dsss93}
\int_{\Omega}|\partial_{x_i}u(x,t)|^{p-1} \sgn(\partial_{x_i}u(x,t))\partial_2f(x, u(x,t)) \partial_{x_i} u(x,t) dx\leqslant (k_fr_\delta+c_f)\|\partial_{x_i}u(\cdot,t)\|_{L^{p}(\Omega)}^{p}.
\end{equation}

Hence, from \eqref{mkqa2da1}-\eqref{dp3dsss93} we conclude that
\[
\frac{d}{dt} \|\partial_{x_i} u(\cdot,t)\|_{L^{p}(\Omega)}^p\leqslant p\|\partial_{x_i}u(\cdot,t)\|_{L^{p}(\Omega)}^{p} \left[-h_0 +k_fr_\delta+c_f  +\frac{M}{\|\partial_{x_i}u(\cdot,t)\|_{L^{p}(\Omega)}}\right],
\]
where $M:=\frac{1}{p}\left[r_\delta\|\partial_{x_i} J\|_{p'}\left(k_g r_\delta + c_g|\Omega|^{\frac{1}{p}} \right)+\left(k_f r_\delta + c_f|\Omega|^{\frac{1}{p}} \right)+pr_\delta\|h\|_{W^{1,\infty}(\Omega)}\right]$.

Let   $\varepsilon =h_0 - (k_fr_\delta+c_f)  >0$. Then, if we consider
\[
\|\partial_{x_i}u(\cdot,t)\|_{L^{p}(\Omega)} \geqslant \dfrac{1}{\varepsilon}M(1+\mu),
 \]
 we have
 \[
\begin{split}
\frac{d}{dt}\|\partial_{x_i}u(\cdot,t)\|_{L^{p}(\Omega)}^{p} &\leqslant
p \|\partial_{x_i}u(t,\cdot)\|_{L^{p}(\Omega)}^{p}\left(-\varepsilon +
\frac{\varepsilon}{1+\mu}\right)\\
&=- \frac{\mu \varepsilon p}{1+\delta} \|\partial_{x_i}u(t,\cdot)\|_{L^{p}(\Omega)}^{p}.
\end{split}
\]
Therefore
\begin{equation}\label{boundsol2}
\begin{split}
\|\partial_{x_i}u(\cdot,t)\|_{L^{p}(\Omega)}^{p} &\leqslant e^{-\frac{ \mu \varepsilon p}{(1+\mu)}t} \|\partial_{x_i}  u_0\|_{L^{p}(\Omega)}^{p} \\
&=e^{- \frac{\mu p}{(1+\mu)}(h_0 -k_fr_\delta-c_f)\,t}\|\partial_{x_i} u_0\|_{L^{p}(\Omega)}^{p}.
\end{split}
\end{equation}
From this and Lemma \ref{L_PullAbs}, we can  conclude that for any $\mu>0$ 
there exists a ball centered at the origin which absorbs bounded subsets of 
$W^{1,p}(\Omega)$  with respect to the nonlinear semigroup $S_J(\cdot)$ 
generated by \eqref{CP}, and therefore the  result follows easily for $1 \leqslant 
p < \infty$ by \cite[Theorem 2.1]{carvalho2012attractors}, and this completes 
the proof of the theorem.
 \qed

\section{Upper semicontinuity of the global attractors}\label{upper-sc}

In this section we will prove the upper semicontinuity of the global attractors $\mathcal{A}_J$ with respect to $J$. For simplicity of notation we denote by $\mathcal{J}$ the class of all $J$ under the conditions in the previous sections.

We know that the nonlinear semigroup in $L^p(\Omega)$, $\{S_J(t);t\geqslant0\}$, generated by \eqref{CP} is given by
\[
S_J(t)u_0(x):=u_J(x,t),
\]
for any $x\in\Omega$ and $t\geqslant0$,  where
\[
u_J(x,t)=e^{-h(x)t}u_0+\int_0^t e^{-h(x)(t-s)}[g(K_J\, u_J(x,s))+f(x,u_J(x,s))]ds.
\]

\begin{theorem}\label{tend}
Under the  conditions of  Theorem \ref{TExistGloAt}. Fixed $J_0\in\mathcal{J}$, for initial data of the Cauchy problem \eqref{CP} in a bounded subset of $L^p(\Omega)$, with $1\leqslant p<\infty$, we have that $u_J$ converges to $u_{J_0}$ in $L^p(\Omega)$ as $J$ converges to $J_0$ in $L^1(\Omega)$.
\end{theorem}

\proof
Note that for any $t\geqslant0$ and $u_0\in B$ with $B\subset L^p(\Omega)$ bounded we have
\[
\begin{split}
\|u_J(\cdot,t)-u_{J_0}(\cdot,t)\|_{L^p(\Omega)}&\leqslant\int_0^te^{-h_0(t-s)} \|g(K_J\, u_J(\cdot,s))-g(K_{J_0} u_{J}(\cdot,s))\|_{L^p(\Omega)}ds\\
&+\int_0^te^{-h_0(t-s)} \|g(K_{J_0} u_J(\cdot,s))-g(K_{J_0} u_{J_0}(\cdot,s))\|_{L^p(\Omega)}ds\\
&+\int_0^te^{-h_0(t-s)} \|f(\cdot,u_J(\cdot,s))-f(\cdot,u_{J_0}(\cdot,s))\|_{L^p(\Omega)}ds.
\end{split}
\]
Now, arguing as in the proof of Proposition \ref{WellP} we obtain a bounded set $D$ which contains $u_J(x,s)$, $u_{J_0}(x,s)$, $K_J\, u_J(x,s)$ and $K_{J_0} u_J(x,s)$. Then using the fact that $g$ and $f$ are Lipschitz continuous on bounded sets, we get
\[
\begin{split}
\|u_J(\cdot,t)-u_{J_0}(\cdot,t)\|_{L^p(\Omega)}&\leqslant L_g\int_0^te^{-h_0(t-s)} \|(K_J\,  - K_{J_0} )u_{J}(\cdot,s)\|_{L^p(\Omega)}ds\\
&+L_g\int_0^te^{-h_0(t-s)} \|K_{J_0} (u_J(\cdot,s)- u_{J_0}(\cdot,s))\|_{L^p(\Omega)}ds\\
&+L_f\int_0^te^{-h_0(t-s)} \| u_J(\cdot,s)- u_{J_0}(\cdot,s)\|_{L^p(\Omega)}ds,
\end{split}
\]
where $L_g$ and $L_f$ are Lipschitz constants of $g$ and $f$, respectively, on $D$. From Lemma \ref{boundK} we have
\[
\begin{split}
\|u_J(\cdot,t)-u_{J_0}(\cdot,t)\|_{L^p(\Omega)}&\leqslant L_g\|J  - J_0\|_{1}\int_0^te^{-h_0(t-s)}  \| u_J(\cdot,s)\|_{L^p(\Omega)}ds\\
&+(L_g+L_f)\int_0^te^{-h_0(t-s)} \| u_J(\cdot,s)- u_{J_0}(\cdot,s)\|_{L^p(\Omega)}ds.
\end{split}
\]
Thanks to \eqref{boundsol}, we have that $u_J(\cdot,s)$ is bounded in $L^p(\Omega)$ and there exists a positive constant $C$ such that
\[
e^{h_0t}\|u_J(\cdot,t)-u_{J_0}(\cdot,t)\|_{L^p(\Omega)}\leqslant \frac{L_g}{h_0}C\|J  - J_0\|_{1}+(L_g+L_f)\int_0^te^{h_0s} \| u_J(\cdot,s)- u_{J_0}(\cdot,s)\|_{L^p(\Omega)}ds
\]
and finally by Gr\"onwall's Lemma we obtain
\[
\|u_J(\cdot,t)-u_{J_0}(\cdot,t)\|_{L^p(\Omega)}\leqslant C_0\|J  - J_0\|_{1}e^{(L_g+L_f-h_0)t},
\]
where $C_0=\frac{L_g}{h_0}C$ for any $t\geqslant0.$
\qed

\begin{remark}\label{dsnr34}
Fixed $J_0\in\mathcal{J}$, for $J$ sufficiently near to $J_0$ in $L^1(\Omega)$, the family of global attractors $\{\mathcal{A}_J; J\in\mathcal{J}\}$ is uniformly bounded in $J$. Indeed, since $\mathcal{A}_J$ is contained in a ball with radius which depend continuously in $J$, we can conclude that there exists a bounded subset of $L^p(\Omega)$ which contains the attractors $\mathcal{A}_J$.
\end{remark}

\begin{theorem}
Under same hypotheses of Theorem  \ref{tend} the family of global attractors $\{\mathcal{A}_J; J\in\mathcal{J}\}$ is upper semicontinuous at $J=J_0$.
\end{theorem}

\proof
Note that, using the invariance of attractors, we have $S_{J}(t)\mathcal{A}_{J}=\mathcal{A}_{J}$ for any $J\in\mathcal{J}$ and
\[
\begin{split}
\operatorname{dist}_{H}(\mathcal{A}_J,\mathcal{A}_{J_0})
&\leqslant \operatorname{dist}_{H}(S_J(t)\mathcal{A}_J,S_{J_0}(t)\mathcal{A}_J)+\operatorname{dist}_{H}(S_{J_0}(t)\mathcal{A}_J,\mathcal{A}_{J_0})\\
&=\sup_{a_J\in \mathcal{A}_J} \operatorname{dist}(S_J(t)a_J,S_{J_0}(t)a_J)+\operatorname{dist}_{H}(S_{J_0}(t)\mathcal{A}_J,\mathcal{A}_{J_0}).
\end{split}
\]
For each $\varepsilon>0$, thanks to Theorem \ref{tend} we have
\[
\sup_{a_J\in \mathcal{A}_J} \operatorname{dist}_H(S_J(t)a_J,S_{J_0}(t)a_J)<\dfrac{\varepsilon}{2},
\]
for any $J$ sufficiently near to $J_0$ in  $L^1$, by the definition of global attractor and Remark \ref{dsnr34}, we have that for some bounded subset $B_0$ of $L^p(\Omega)$ we obtain 
\[
\operatorname{dist}_{H}(S_{J_0}(t)\mathcal{A}_J,\mathcal{A}_{J_0})\leqslant \operatorname{dist}_{H}(S_{J_0}(t) B_0,\mathcal{A}_{J_0})<\dfrac{\varepsilon}{2},
\]
for any $t$ sufficiently large. Therefore, for  $J$ sufficiently near to $J_0$ in  $L^1$, and for any $t$ sufficiently large, we get
\[
\operatorname{dist}_{H}(\mathcal{A}_J,\mathcal{A}_{J_0})<\varepsilon.
\]
\qed

\bibliographystyle{abbrv}

\begin{thebibliography}{10}

\bibitem{amari1977dynamics}
S.-i. Amari.
\newblock Dynamics of pattern formation in lateral-inhibition type neural
  fields.
\newblock {\em Biological cybernetics}, 27(2):77--87, 1977.

\bibitem{Mazon_Rossi_2008}
F.~Andreu, J.~M. Maz\'on, J.~D. Rossi, and J.~Toledo.
\newblock The {N}eumann problem for nonlocal nonlinear diffusion equations.
\newblock {\em J. Evol. Equ.}, 8(1):189--215, 2008.

\bibitem{Rossi_libro}
F.~Andreu-Vaillo, J.~M. Maz\'on, J.~D. Rossi, and J.~J. Toledo-Melero.
\newblock {\em Nonlocal diffusion problems}, volume 165 of {\em Mathematical
  Surveys and Monographs}.
\newblock American Mathematical Society, Providence, RI; Real Sociedad
  Matem\'atica Espa\~nola, Madrid, 2010.

\bibitem{Berestycki}
H.~Berestycki, G.~Nadin, B.~Perthame, and L.~Ryzhik.
\newblock The non-local {F}isher-{KPP} equation: travelling waves and steady
  states.
\newblock {\em Nonlinearity}, 22(12):2813--2844, 2009.

\bibitem{bezerra2012existence}
F.~D. Bezerra, A.~L. Pereira, and S.~H. Da~Silva.
\newblock Existence and continuity of global attractors and nonhomogeneous
  equilibria for a class of evolution equations with non local terms.
\newblock {\em Journal of Mathematical Analysis and Applications},
  396(2):590--600, 2012.

\bibitem{carvalho2012attractors}
A.~Carvalho, J.~A. Langa, and J.~Robinson.
\newblock {\em Attractors for infinite-dimensional non-autonomous dynamical
  systems}, volume 182.
\newblock Springer Science \& Business Media, 2012.

\bibitem{chasseigne2006asymptotic}
E.~Chasseigne, M.~Chaves, and J.~D. Rossi.
\newblock Asymptotic behavior for nonlocal diffusion equations.
\newblock {\em Journal de math{\'e}matiques pures et appliqu{\'e}es},
  86(3):271--291, 2006.

\bibitem{chasseigne2013nonlocal}
E.~Chasseigne, S.~Sastre-Gomez, et~al.
\newblock A nonlocal two phase stefan problem.
\newblock {\em Differential and Integral Equations}, 26(11/12):1335--1360,
  2013.

\bibitem{cortazar2005nonlocal}
C.~Cortazar, M.~Elgueta, and J.~D. Rossi.
\newblock A nonlocal diffusion equation whose solutions develop a free
  boundary.
\newblock In {\em Annales Henri Poincar{\'e}}, volume~6, pages 269--281.
  Springer, 2005.

\bibitem{cortazar2008approximate}
C.~Cortazar, M.~Elgueta, J.~D. Rossi, and N.~Wolanski.
\newblock How to approximate the heat equation with neumann boundary conditions
  by nonlocal diffusion problems.
\newblock {\em Archive for Rational Mechanics and Analysis}, 187(1):137--156,
  2008.

\bibitem{da2012properties}
S.~H. Da~Silva.
\newblock Properties of an equation for neural fields in a bounded domain.
\newblock {\em Electronic Journal of Differential Equations}, 2012(42):1--9,
  2012.

\bibitem{da2013finite}
S.~H. Da~Silva and F.~D. Bezerra.
\newblock Finite fractal dimensionality of attractors for nonlocal evolution
  equations.
\newblock {\em Electronic Journal of Differential Equations}, 2013(221):1--9,
  2013.

\bibitem{da2009global}
S.~H. Da~Silva and A.~L. Pereira.
\newblock Global attractors for neural fields in a weighted space.
\newblock {\em Matem{\'a}tica Contemporanea}, 36:139--153, 2009.

\bibitem{daleckii2002stability}
J.~L. Daleckii and M.~G. Kre\_n.
\newblock {\em Stability of solutions of differential equations in Banach
  space}.
\newblock Number~43. American Mathematical Soc., 2002.

\bibitem{de1995travelling}
A.~De~Masi, T.~Gobron, and E.~Presutti.
\newblock Travelling fronts in non-local evolution equations.
\newblock {\em Archive for rational mechanics and analysis}, 132(2):143--205,
  1995.

\bibitem{de2000critical}
A.~De~Masi, E.~Olivieri, and E.~Presutti.
\newblock Critical droplet for a non local mean field equation.
\newblock {\em Markov Processes and Related Fields}, 6:439--472, 2000.

\bibitem{de1994glauber}
A.~De~Masi, E.~Orlandi, E.~Presutti, and L.~Triolo.
\newblock Glauber evolution with kac potentials. i. mesoscopic and macroscopic
  limits, interface dynamics.
\newblock {\em Nonlinearity}, 7(3):633, 1994.

\bibitem{de1994stability}
A.~De~Masi, E.~Orlandi, E.~Presutti, and L.~Triolo.
\newblock Stability of the interface in a model of phase separation.
\newblock {\em Proceedings of the Royal Society of Edinburgh Section A:
  Mathematics}, 124(5):1013--1022, 1994.

\bibitem{de1994uniqueness}
A.~De~Masi, E.~Orlandi, E.~Presutti, and L.~Triolo.
\newblock Uniqueness of the instanton profile and global stability in non local
  evolution equations.
\newblock {\em Rendiconti di Matematica e Delle sue Applicazioni}, 4:693--723,
  1994.

\bibitem{folland1995introduction}
G.~B. Folland.
\newblock {\em Introduction to partial differential equations}.
\newblock Princeton university press, 1995.

\bibitem{hutson}
V.~Hutson, S.~Martinez, K.~Mischaikow, and G.~T. Vickers.
\newblock The evolution of dispersal.
\newblock {\em J. Math. Biol.}, 47(6):483--517, 2003.

\bibitem{ladas1972differential}
G.~E. Ladas and V.~Lakshmikantham.
\newblock {\em Differential equations in abstract spaces}.
\newblock Elsevier, 1972.

\bibitem{pereira2006global}
A.~L. Pereira.
\newblock Global attractor and nonhomogeneous equilibria for a nonlocal
  evolution equation in an unbounded domain.
\newblock {\em Journal of Differential Equations}, 226(1):352--372, 2006.

\bibitem{rall2014nonlinear}
L.~B. Rall.
\newblock {\em Nonlinear functional analysis and applications: proceedings of
  an advanced seminar conducted by the Mathematics Research Center, the
  University of Wisconsin, Madison, October 12-14, 1970}.
\newblock Number~26. Elsevier, 2014.

\bibitem{rodriguez2014nonlinear}
A.~Rodr{\'\i}guez-Bernal and S.~Sastre-G{\'o}mez.
\newblock Nonlinear nonlocal reaction-diffusion equations.
\newblock In {\em Advances in Differential Equations and Applications}, pages
  53--61. Springer, 2014.

\bibitem{rodriguez-bernal_sastre-gomez_2016}
A.~Rodr{\'i}guez-Bernal and S.~Sastre-Gomez.
\newblock Linear non-local diffusion problems in metric measure spaces.
\newblock {\em Proceedings of the Royal Society of Edinburgh: Section A
  Mathematics}, 146(4):833?863, 2016.

\bibitem{ref6}
E.~Valdinoci.
\newblock From the long jump random walk to the fractional {L}aplacian.
\newblock {\em Bol. Soc. Esp. Mat. Apl. S$\vec{\rm e}$MA}, 49:33--44, 2009.

\end{thebibliography}

\end{document}